\documentclass[reqno,b5paper]{amsart}

\usepackage{amsmath}
\usepackage{amssymb}
\usepackage{amsthm}
\usepackage{enumerate}
\usepackage[mathscr]{eucal}
\usepackage{eqlist}

\theoremstyle{plain}
\newtheorem{theorem}{Theorem}[section]

\newtheorem{claim}{Claim}

\theoremstyle{definition}

\numberwithin{equation}{section}

\numberwithin{equation}{section}

\begin{document}
\title[Non-linear new product $A^{*}B - B^{*}A$ derivations on $\ast$-algebras]%
{Non-linear new product $A^{*}B - B^{*}A$ derivations on $\ast$-algebras}
\author[A. Taghavi, M.Razeghi]%
{Ali Taghavi, Mehran Razeghi}

\newcommand{\acr}{\newline\indent}
\address{\llap{*\,}Department of Mathematics\\ Faculty of Mathematical
Sciences\\ University of Mazandaran\\ P. O. Box 47416-1468\\
Babolsar, Iran.} 
 \email{taghavi@umz.ac.ir, razeghi.mehran19@yahoo.com}


\begin{abstract}
Let $\mathcal{A}$ be a prime $\ast$-algebra. In this paper, we suppose that $\Phi:\mathcal{A}\to\mathcal{A}$ satisfies
$$\Phi(A\diamond B)=\Phi(A)\diamond B+A\diamond\Phi(B)$$
where  $A\diamond B = A^{*}B - B^{*}A$ for all $A,B\in\mathcal{A}$ .We will show that if $\Phi(\alpha \frac{I}{2})$ is self-adjoint for $\alpha\in\{1,i\}$ then $\Phi$ is additive $\ast$-derivation.   
\end{abstract}

\maketitle

\section{Introduction}\label{intro}

Let $\mathcal{R}$ be a $*$-algebra. For $A,B\in\mathcal{R}$, denoted by
$A\bullet B=AB+BA^{*}$ and $[A,B]_{*}=AB-BA^{*}$, which are $\ast$-Jordan product and $\ast$-Lie product, respectively. These products are found playing a
more and more important role in some research topics, and its study
has recently attracted many author's attention (for example, see
\cite{cui,li2,mol,taghavi}).

Recall that a map $\Phi:\mathcal{R}\to\mathcal{R}$ is said to be an additive derivation if $$\Phi(A+B)=\Phi(A)+\Phi(B)$$ and $$\Phi(AB)=\Phi(A)B+A\Phi(B)$$
for all $A,B\in\mathcal{R}$. A map $\Phi$ is additive $\ast$-derivation if it is an additive derivation and $\Phi(A^{*})=\Phi(A)^{*}$. 
Derivations are very important maps both in theory and applications, and have been studied intensively (\cite{chr,sak,sem1,sem2}).

Let us define $\lambda$-Jordan $\ast$-product by $A\bullet_{\lambda}B =
AB + \lambda BA^{\ast}$. We say that the map $\Phi$ with the property of
$\Phi(A\bullet_{\lambda} B) = \Phi(A)\bullet_{\lambda}B +  A
\bullet_{\lambda}\Phi(B)$ is a $\lambda$-Jordan $\ast$-derivation
map. It is clear that for $\lambda = - 1$ and $\lambda = 1,$ the $\lambda
$-Jordan $\ast$-derivation map is a $\ast$-Lie derivation and
$\ast$-Jordan derivation, respectively \cite{bai}. 

A von Neumann algebra $\mathcal{A}$ is a self-adjoint subalgebra of some $B(H)$, the algebra of bounded linear operators acting on a complex Hilbert space, which satisfies the double commutant property: $\mathcal{A}^{''}=\mathcal{A}$ where $\mathcal{A}^{'}=\{T\in B(H), TA=AT, \forall A\in\mathcal{A}\}$ and $\mathcal{A}^{''}=\{\mathcal{A}^{'}\}^{'}$. Denote by $\mathcal{Z}(\mathcal{A})=\mathcal{A}^{'}\cap \mathcal{A}$ the center of $\mathcal{A}$. A von Neumann algebra $\mathcal{A}$ is called a factor if its center is trivial, that is, $\mathcal{Z}(\mathcal{A})=\mathbb{C}I$. For $A\in\mathcal{A}$, recall that the central carrier of $A$, denoted by $\overline{A}$, is the smallest central projection $P$ such that $PA=A$. It is not difficult to see that $\overline{A}$ is the projection onto the closed subspace spanned by $\{BAx : B\in \mathcal{A}, x\in H\}$. If $A$ is self-adjoint, then the core of $A$, denoted by $\underline{A}$, is $\sup\{S\in\mathcal{Z}(\mathcal{A}): S=S^{*}, S\leq A\}$. If $A=P$ is a projection, it is clear that $\underline{P}$ is the largest central projection $Q$ satisfying $Q\leq P$. A projection $P$ is said to be core-free if $\underline{P}=0$ (see \cite{mie2}). It is easy to see that $\underline{P}=0$ if and only if $\overline{I-P}=I$, \cite{kad1,kad2}.

Recently, Yu and Zhang in \cite{yu}
 proved that every non-linear $\ast$-Lie derivation  from a factor von Neumann algebra into itself is an additive $\ast$-derivation. Also,
 Li, Lu and Fang
in \cite{li} have investigated a non-linear $\lambda$-Jordan
$\ast$-derivation. They showed that if
$\mathcal{A}\subseteq\mathcal{B(H)}$ is a von Neumann algebra
without central abelian projections and $\lambda$ is a non-zero scalar,
then $\Phi:\mathcal{A} \longrightarrow \mathcal{B(H)}$ is a
non-linear $\lambda$-Jordan $\ast$-derivation if and only if $\Phi$ is
an additive $\ast$-derivation.

On the other hand, many mathematician devoted themselves to study the $\ast$-Jordan product $A\bullet B=AB+BA^{*}$. In \cite{zha}, F. Zhang proved that every non-linear $\ast$-Jordan derivation map $\Phi:\mathcal{A}\to\mathcal{A}$ on a factor von neumann algebra with $I_{\mathcal{A}}$ the identity of it is an additive $\ast$-derivation.

 In \cite{taghavi2}, we showed that $\ast$-Jordan derivation map on every
factor von Neumann algebra $\mathcal{A}
 \subseteq \mathcal{B(H)}$  is
additive $\ast$-derivation.

Very recently the authors of \cite{chen} discussed some bijective maps preserving the new product $A^{*}B+B^{*}A$ between von Neumann algebras with no central abelian projections. In other words, $\Phi$ holds in the following condition
$$\Phi(A^{*}B+B^{*}A)=\Phi(A)^{*}\Phi(B)+\Phi(B)^{*}\Phi(A).$$
 They showed that such a map is sum of a linear $\ast$-isomorphism and a conjugate linear $\ast$-isomorphism. 

Motivated by the above results, in this paper, we prove that if $\mathcal{A}$ is a prime $\ast$-algebra then $\Phi:\mathcal{A}\to\mathcal{A}$ which holds in the following condition
$$\Phi(A\diamond B)=\Phi(A)\diamond B+A\diamond\Phi(B)$$
where  $A\diamond B = A^{*}B - B^{*}A$ for all $A,B\in\mathcal{A}$, is additive $\ast$-derivation.

We say that $\mathcal{A}$ is prime, that is, for $A,B \in \mathcal{A}$ if $A\mathcal{A}B =\lbrace0\rbrace,$ then $A = 0$ or $B = 0$.

\section{Main Results}

Our main theorem is as follows:
\begin{theorem}\label{mainsh}
Let $\mathcal{A}$ be a prime $\ast$-algebra. Let $\Phi:\mathcal{A}\to \mathcal{A}$ satisfies in  
\begin{equation}\label{sh1}
\Phi(A\diamond B)=\Phi(A)\diamond B+A\diamond\Phi(B) 
\end{equation}
where  $A\diamond B = A^{*}B - B^{*}A$ for all $A,B\in\mathcal{A}$. If $\Phi(\alpha \frac{I}{2})$ is self-adjoint operator for $\alpha\in\{1,i\}$ then $\Phi$ is additive $\ast$-derivation.   
\end{theorem}
\textbf{Proof.}
Let $P_{1}$ be a nontrivial
projection in $\mathcal{A}$ and $P_{2}=I_{\mathcal{A}}-P_{1}$.
Denote $\mathcal{A}_{ij}=P_{i}\mathcal{A}P_{j},\ i,j=1,2,$ then
$\mathcal{A}=\sum_{i,j=1}^{2}\mathcal{A}_{ij}$. For every
$A\in\mathcal{A}$ we may write $A=A_{11}+A_{12}+A_{21}+A_{22}$. In
all that follow, when we write $A_{ij}$, it indicates that
$A_{ij}\in\mathcal{A}_{ij}$. For showing additivity of $\Phi$ on $\mathcal{A}$, we  use above
partition of $\mathcal{A}$ and give some claims that prove $\Phi$ is
additive on each $\mathcal{A}_{ij}, \ i,j=1,2$.\\

We prove the above theorem by several claims.
\begin{claim}\label{cl1}
We show that $\Phi(0)=0$.
\end{claim}
This claim is easy to prove.

\begin{claim}\label{cl103}
$\Phi(i\frac{I}{2})=\Phi(\frac{I}{2})=\Phi(-i\frac{I}{2})=0$.
\end{claim}
Consider $\Phi(\frac{I}{2} \diamond i\frac{I}{2})=\Phi(\frac{I}{2}) \diamond i\frac{I}{2}+\frac{I}{2} \diamond \Phi(i \frac{I}{2})$
that imply
\begin{equation}\label{l4}
\Phi(i \frac{I}{2})=\frac{i}{2}\Phi(\frac{I}{2})^{*}+\frac{i}{2} \Phi(\frac{I}{2})+\frac{1}{2} \Phi(i \frac{I}{2})-\frac{1}{2} \Phi(i \frac{I}{2})^{*}=i\Phi(\frac{I}{2}).
\end{equation}
By taking the adjoint of above equation we have $\Phi(i\frac{I}{2})=\Phi(\frac{I}{2})=0$
\\Consider $\Phi(i\frac{I}{2} \diamond \frac{I}{2})=\Phi(i\frac{I}{2}) \diamond \frac{I}{2}+i\frac{I}{2} \diamond \Phi( \frac{I}{2})$
that imply
\begin{equation}\label{z1}
\Phi(-i \frac{I}{2})=\frac{I}{2}\Phi(i\frac{I}{2})-\frac{I}{2} \Phi(i\frac{I}{2})-\frac{i}{2} \Phi( \frac{I}{2})-\frac{i}{2} \Phi( \frac{I}{2})^{*}.
\end{equation}
Since $\Phi(\frac{I}{2})=\Phi(i \frac{I}{2})=0$, so we have $\Phi(-i\frac{I}{2})=0$.
\begin{claim}
$\Phi(-\frac{I}{2})=0$.
\end{claim}
Consider $\Phi(-\frac{I}{2} \diamond i\frac{I}{2})=\Phi(-\frac{I}{2}) \diamond i\frac{I}{2}$. So, we have
$$\Phi(-\frac{I}{2})^{*}+ \Phi(-\frac{I}{2})=0.$$
It follows that
\begin{equation}\label{z2}
\Phi(-\frac{I}{2})^{*}=- \Phi(-\frac{I}{2}).
\end{equation}
Also
$$\Phi(\frac{I}{2} \diamond - \frac{I}{2})=\frac{I}{2} \diamond \Phi(-\frac{I}{2}).$$
So $\Phi(- \frac{I}{2})- \Phi(- \frac{I}{2})^{*}=0$, then we have
\begin{equation}\label{z3}
\Phi(- \frac{I}{2})= \Phi(- \frac{I}{2})^{*}
\end{equation}
from (\ref{z2}), (\ref{z3}) we have $\Phi(-\frac{I}{2})=0.$
\begin{claim}\label{cl3}
For each $A \in \mathcal{A}$, we have
\begin{enumerate}
\item
$\Phi(-iA)=-i\Phi(A).$
\item
$\Phi(iA)=i\Phi(A).$
\end{enumerate}
\end{claim}
We can check to see that 
$$\Phi\left(-iA\diamond \frac{I}{2}\right)=\Phi\left(A\diamond i\frac{I}{2}\right).$$
So, 
$$\Phi(-iA)\diamond\frac{I}{2}=\Phi(A)\diamond i\frac{I}{2}.$$
It follows that
\begin{equation}\label{eq1cl3}
\Phi(-iA)^{*}-\Phi(-iA)=i\Phi(A)^{*}+i\Phi(A).
\end{equation}
On the other hand, one can check that
$$\Phi\left(-iA\diamond i\frac{I}{2}\right)=\Phi\left(A \diamond -\frac{I}{2} \right).$$
So, $$\Phi(-iA)\diamond i\frac{I}{2}=\Phi(A)\diamond -\frac{I}{2}.$$
It follows that
\begin{equation}\label{eq2cl3}
i\Phi(-iA)^{*}+i\Phi(-iA)=-\Phi(A)^{*}+\Phi(A).
\end{equation}
Equivalently, we obtain
\begin{equation}\label{eq3cl3}
-\Phi(-iA)^{*}-\Phi(-iA)=-i\Phi(A)^{*}+i\Phi(A).
\end{equation}
By adding equations (\ref{eq1cl3}) and (\ref{eq3cl3}) we have
$$\Phi(-iA)=-i\Phi(A).$$

Similarly, we can show that $\Phi(iA)=i\Phi(A)$.

\begin{claim}\label{cl4}
For each $A_{11} \in \mathcal{A}_{11}$, $A_{12} \in \mathcal{A}_{12}$ we have 
$$\Phi(A_{11}+A_{12})=\Phi(A_{11})+\Phi(A_{12}).$$
\end{claim}
Let $T=\Phi(A_{11}+A_{12})-\Phi(A_{11})-\Phi(A_{12})$, we should prove that $T=0$.\\
For $X_{21}\in\mathcal{A}_{21}$ we can write that
\begin{eqnarray*}
&&\Phi(A_{11}+A_{12})\diamond X_{21}+(A_{11}+A_{12})\diamond\Phi(X_{21})=\Phi((A_{11}+A_{12})\diamond X_{21})\\
&&=\Phi(A_{11}\diamond X_{21})+\Phi(A_{12}\diamond X_{21})=\Phi(A_{11})\diamond X_{21}+A_{11}\diamond\Phi(X_{21})\\
&&+\Phi(A_{12})\diamond X_{21}+A_{12}\diamond\Phi(X_{21})\\
&&=(\Phi(A_{11})+\Phi(A_{12}))\diamond X_{21}+(A_{11}+A_{12})\diamond\Phi(X_{21}).
\end{eqnarray*}
So, we obtain 
$$T\diamond X_{21}=0.$$
Since $T=T_{11}+T_{12}+T_{21}+T_{22}$ we have
$$T_{21}^{*}X_{21}+T_{22}^{*}X_{21}-X_{21}^{*}T_{21}-X_{21}^{*}T_{22}=0.$$
From the above equation and primeness of $\mathcal{A}$ we have
$T_{22}=0$ and 
\begin{equation}\label{eq4cl3}
T_{21}^{*}X_{21}-X_{21}^{*}T_{21}=0.
\end{equation}
On the other hand, similarly by applying $iX_{21}$ instead of $X_{21}$ in above, we obtain 
$$iT_{21}^{*}X_{21}+iT_{22}^{*}X_{21}+iX_{21}^{*}T_{21}+iX_{21}^{*}T_{22}=0.$$
Since $T_{22}=0$ we obtain from the above equation that
\begin{equation}\label{eq5cl3}
-T_{21}^{*}X_{21}-X_{21}^{*}T_{21}=0.
\end{equation}
From (\ref{eq4cl3}) and (\ref{eq5cl3}) we have
$$X_{21}^{*}T_{21}=0.$$
Since $\mathcal{A}$ is prime, then we get $T_{21}=0$.\\

It suffices to show that $T_{12}=T_{11}=0$. For this purpose for $X_{12}\in\mathcal{A}_{12}$ we write
\begin{eqnarray*}
&&\Phi(((A_{11}+A_{12})\diamond X_{12})\diamond P_{1})=\Phi((A_{11}+A_{12})\diamond X_{12})\diamond P_{1}+((A_{11}+A_{12})\diamond X_{12})\diamond \Phi(P_{1})\\
&&=(\Phi(A_{11}+A_{12})\diamond X_{12}+(A_{11}+A_{12})\diamond\Phi(X_{12}))\diamond P_{1}+(A_{11}+A_{12})\diamond X_{12}\diamond\Phi(P_{1})\\
&&=\Phi(A_{11}+A_{12})\diamond X_{12}\diamond P_{1}+A_{11}\diamond\Phi(X_{12})\diamond P_{1}+A_{12}\diamond\Phi(X_{12})\diamond P_{1}\\
&&+A_{11}\diamond X_{12}\diamond\Phi(P_{1})+A_{12}\diamond X_{12}\diamond \Phi(P_{1}).
\end{eqnarray*}
So, we showed that 
\begin{eqnarray}
&&\Phi(((A_{11}+A_{12})\diamond X_{12})\diamond P_{1})=\Phi(A_{11}+A_{12})\diamond X_{12}\diamond P_{1}+A_{11}\diamond\Phi(X_{12})\diamond P_{1}\nonumber\\&&+A_{12}\diamond \Phi(X_{12})\diamond P_{1}+A_{11}\diamond X_{12}\diamond \Phi(P_{1})+A_{12}\diamond X_{12}\diamond \Phi(P_{1}). \label{eq6cl3}
\end{eqnarray}
Since $A_{12}\diamond X_{12}\diamond P_{1}=0$ we have
\begin{eqnarray*}
&&\Phi(((A_{11}+A_{12})\diamond X_{12})\diamond P_{1})=\Phi((A_{11}\diamond X_{12})\diamond P_{1})+\Phi((A_{12}\diamond X_{12})\diamond P_{1})\\
&&=\Phi(A_{11}\diamond X_{12})\diamond P_{1}+(A_{11}\diamond X_{12})\diamond \Phi(P_{1})+\Phi(A_{12}\diamond X_{12})\diamond P_{1}+(A_{12}\diamond X_{12})\diamond \Phi(P_{1})\\
&&=(\Phi(A_{11})\diamond X_{12}+A_{11}\diamond\Phi(X_{12}))\diamond P_{1}+(A_{11}\diamond X_{12})\diamond \Phi(P_{1})\\
&&+(\Phi(A_{12})\diamond X_{12}+A_{12}\diamond \Phi(X_{12}))\diamond P_{1}+(A_{12}\diamond X_{12})\diamond \Phi(P_{1})\\
&&=\Phi(A_{11})\diamond X_{12}\diamond P_{1}+A_{11}\diamond \Phi(X_{12})\diamond P_{1}+A_{11}\diamond X_{12}\diamond\Phi(P_{1})\\
&&+\Phi(A_{12})\diamond X_{12}\diamond P_{1}+A_{12}\diamond \Phi(X_{12})\diamond P_{1}+A_{12}\diamond X_{12}\diamond \Phi(P_{1}).
\end{eqnarray*}
So, 
\begin{eqnarray}
&&\Phi(((A_{11}+A_{12})\diamond X_{12})\diamond P_{1})=\Phi(A_{11})\diamond X_{12}\diamond P_{1}+A_{11}\diamond \Phi(X_{12})\diamond P_{1}\nonumber\\
&&+A_{11}\diamond X_{12}\diamond\Phi(P_{1})+\Phi(A_{12})\diamond X_{12}\diamond P_{1}\nonumber\\
&&+A_{12}\diamond \Phi(X_{12})\diamond P_{1}+A_{12}\diamond X_{12}\diamond \Phi(P_{1}). \label{eq7cl3}
\end{eqnarray}
From (\ref{eq6cl3}) and (\ref{eq7cl3}) we have
$$\Phi(A_{11}+A_{12})\diamond X_{12}\diamond P_{1}=\Phi(A_{11})\diamond X_{12}\diamond P_{1}+\Phi(A_{12})\diamond X_{12}\diamond P_{1}.$$
It follows that $T\diamond X_{12}\diamond P_{1}=0$, so $T_{11}^{*}X_{12}-X_{12}^{*}T_{11}=0$. We have $T_{11}^*{}X_{12}=0$ or $T_{11}XP_{2}=0$ for all $X\in\mathcal{A}$, then we have
$T_{11}=0$. Similarly, we can show that $T_{12}=0$ by applying $P_{2}$ instead of $P_{1}$ in above.

\begin{claim}\label{cl5}
For each $A_{11}\in\mathcal{A}_{11}, A_{12}\in\mathcal{A}_{12}, A_{21}\in\mathcal{A}_{21}$ and $A_{22}\in\mathcal{A}_{22}$ we have
\begin{enumerate}
\item
$$\Phi(A_{11}+A_{12}+A_{21})=\Phi(A_{11})+\Phi(A_{12})+\Phi(A_{21}).$$
\item
$$\Phi(A_{12}+A_{21}+A_{22})=\Phi(A_{12})+\Phi(A_{21})+\Phi(A_{22}).$$
\end{enumerate}
\end{claim}
We show that 
$$T=\Phi(A_{11}+A_{12}+A_{21})-\Phi(A_{11})-\Phi(A_{12})-\Phi(A_{21})=0.$$
So, we have
\begin{eqnarray*}
&&\Phi(A_{11}+A_{12}+A_{21})\diamond X_{21}+(A_{11}+A_{12}+A_{21})\diamond \Phi(X_{21})\\
&&=\Phi((A_{11}+A_{12}+A_{21})\diamond X_{21})=\Phi(A_{11}\diamond X_{21})+\Phi(A_{12}\diamond X_{21})+\Phi(A_{21}\diamond X_{21})\\
&&=(\Phi(A_{11})+\Phi(A_{12})+\Phi(A_{21}))\diamond X_{21}+(A_{11}+A_{12}+A_{21})\diamond \Phi(X_{21}).
\end{eqnarray*}
It follows that $T\diamond X_{21}=0$. Since $T=T_{11}+T_{12}+T_{21}+T_{22}$ we have 
$$T_{22}^{*}X_{21}+T_{21}^{*}X_{21}-X_{21}^{*}T_{22}-C_{21}^{*}T_{21}=0.$$
Therefore, $T_{22}=T_{21}=0$.\\
From Claim \ref{cl4}, we obtain
\begin{eqnarray*}
&&\Phi(A_{11}+A_{12}+A_{21})\diamond X_{12}+(A_{11}+A_{12}+A_{21})\diamond \Phi(X_{12})\\
&&=\Phi((A_{11}+A_{12}+A_{21})\diamond X_{12})=\Phi((A_{11}+A_{12})\diamond X_{12})+\Phi(A_{21}\diamond X_{12})\\
&&=\Phi(A_{11}\diamond X_{12})+\Phi(A_{12}\diamond X_{12})+\Phi(A_{21}\diamond X_{12})\\
&&=(\Phi(A_{11})+\Phi(A_{12})+\Phi(A_{21}))\diamond X_{12}+(A_{11}+A_{12}+A_{21})\diamond \Phi(X_{12}).
\end{eqnarray*}
Hence, 
$$T_{11}^{*}X_{12}+T_{12}^{*}X_{12}-X_{12}^{*}T_{11}-X_{12}^{*}T_{12}=0.$$
Then $T_{11}=T_{12}=0$.
Similarly  $$\Phi(A_{12}+A_{21}+A_{22})=\Phi(A_{12})+\Phi(A_{21})+\Phi(A_{22}).$$

\begin{claim}\label{cl6}
For each $A_{11}\in\mathcal{A}_{11}, A_{12}\in\mathcal{A}_{12}, A_{21}\in\mathcal{A}_{21}$ and $A_{22}\in\mathcal{A}_{22}$ we have
$$\Phi(A_{11}+A_{12}+A_{21}+A_{22})=\Phi(A_{11})+\Phi(A_{12})+\Phi(A_{21})+\Phi(A_{22}).$$
\end{claim}
We show that 
$$T=\Phi(A_{11}+A_{12}+A_{21}+A_{22})-\Phi(A_{11})-\Phi(A_{12})-\Phi(A_{21})-\Phi(A_{22})=0.$$
From Claim \ref{cl5}, we have
\begin{eqnarray*}
&&\Phi(A_{11}+A_{12}+A_{21}+A_{22})\diamond X_{12}+(A_{11}+A_{12}+A_{21}+A_{22})\diamond\Phi(X_{12})\\
&&=\Phi((A_{11}+A_{12}+A_{21}+A_{22})\diamond X_{12})\\
&&=\Phi((A_{11}+A_{12}+A_{21})\diamond X_{12})+\Phi(A_{22}\diamond X_{12})\\
&&=\Phi(A_{11}\diamond X_{12})+\Phi(A_{12}\diamond X_{12})+\Phi(A_{21}\diamond X_{12})+\Phi(A_{22}\diamond X_{12})\\
&&=(\Phi(A_{11})+\Phi(A_{12})+\Phi(A_{21})+\Phi(A_{22}))\diamond X_{12}\\
&&+(A_{11}+A_{12}+A_{21}+A_{22})\diamond \Phi(X_{12}).
\end{eqnarray*}
So, $T\diamond X_{12}=0$. It follows that
$$T_{11}^{*}X_{12}+T_{12}^{*}X_{12}-X_{12}^{*}T_{11}-X_{12}^{*}T_{12}=0.$$ Then $T_{11}=T_{12}=0$.

Similarly, by applying $X_{21}$ instead of $X_{12}$ in above, we obtain $T_{21}=T_{22}=0$.

\begin{claim}\label{cl7}
For each $A_{ij},B_{ij} \in \mathcal{A}_{ij}$ such that $i\neq j$, we have 
$$\Phi(A_{ij}+B_{ij})=\Phi(A_{ij})+\Phi(B_{ij}).$$
\end{claim}
It is easy to show that
$$(P_{i}+A_{ij})(P_{j}+B_{ij})-(P_{j}+B_{ij}^{*})(P_{i}+A_{ij}^{*})=A_{ij}+B_{ij}-A_{ij}^{*}-B_{ij}^{*}.$$
So, we can write
\begin{eqnarray*}
&&\Phi(A_{ij}+B_{ij})+\Phi(-A_{ij}^{*}-B_{ij}^{*})=\Phi((P_{i}+A_{ij}^{*})\diamond (P_{j}+B_{ij}))\\
&&=\Phi(P_{i}+A_{ij}^{*})\diamond (P_{j}+B_{ij})+(P_{i}+A_{ij}^{*})\diamond \Phi(P_{j}+B_{ij})\\
&&=(\Phi(P_{i})+\Phi(A_{ij}^{*}))\diamond (P_{j}+B_{ij})+(P_{i}+A_{ij}^{*})\diamond (\Phi(P_{j})+\Phi(B_{ij}))\\
&&=\Phi(P_{i})\diamond B_{ij}+P_{i}\diamond \Phi(B_{ij})+\Phi(A_{ij}^{*})\diamond P_{j}+A_{ij}^{*}\diamond \Phi(P_{j})\\
&&=\Phi(P_{i}\diamond B_{ij})+\Phi(A_{ij}^{*}\diamond P_{j})\\
&&=\Phi(B_{ij})+\Phi(-B_{ij}^{*})+\Phi(A_{ij})+\Phi(-A_{ij}^{*}).
\end{eqnarray*}
Therefore, we show that
\begin{equation}\label{eq1cl7}
\Phi(A_{ij}+B_{ij})+\Phi(-A_{ij}^{*}-B_{ij}^{*})=\Phi(A_{ij})+\Phi(B_{ij})+\Phi(-A_{ij}^{*})+\Phi(-B_{ij}^{*}).
\end{equation}
By an easy computation, we can write
$$(P_{i}+A_{ij})(iP_{j}+iB_{ij})-(-iP_{j}-iB_{ij}^{*})(P_{i}+A_{ij}^{*})=iA_{ij}+iB_{ij}+iA_{ij}^{*}+iB_{ij}^{*}.$$
Then, we have
\begin{eqnarray*}
&&\Phi(iA_{ij}+iB_{ij})+\Phi(iA_{ij}^{*}+iB_{ij}^{*})=\Phi((P_{i}+A_{ij}^{*})\diamond (iP_{j}+iB_{ij}))\\
&&=\Phi(P_{i}+A_{ij}^{*})\diamond (iP_{j}+iB_{ij})+(P_{i}+A_{ij}^{*})\diamond \Phi(iP_{j}+iB_{ij})\\
&&=(\Phi(P_{i})+\Phi(A_{ij}^{*}))\diamond(iP_{j}+iB_{ij})+(P_{i}+A_{ij}^{*})(\Phi(iP_{j})+\Phi(iB_{ij}))\\
&&=\Phi(P_{i})\diamond iB_{ij}+P_{i}\diamond \Phi(iB_{ij})+\Phi(A_{ij}^{*})\diamond iP_{j}+A_{ij}^{*}\diamond \Phi(iP_{j})\\
&&=\Phi(P_{i}\diamond iB_{ij})+\Phi(A_{ij}^{*}\diamond iP_{j})\\
&&=\Phi(iB_{ij})+\Phi(iB_{ij}^{*})+\Phi(iA_{ij})+\Phi(iA_{ij}^{*}).
\end{eqnarray*}
We showed that 
$$\Phi(iA_{ij}+iB_{ij})+\Phi(iA_{ij}^{*}+iB_{ij}^{*})=\Phi(iB_{ij})+\Phi(iB_{ij}^{*})+\Phi(iA_{ij})+\Phi(iA_{ij}^{*}).$$
From Claim \ref{cl3} and the above equation, we have
\begin{equation}\label{eq2cl7}
\Phi(A_{ij}+B_{ij})-\Phi(-A_{ij}^{*}-B_{ij}^{*})=\Phi(B_{ij})-\Phi(-B_{ij}^{*})+\Phi(A_{ij})-\Phi(-A_{ij}^{*}).
\end{equation}
By adding equations (\ref{eq1cl7}) and (\ref{eq2cl7}), we obtain
$$\Phi(A_{ij}+B_{ij})=\Phi(A_{ij})+\Phi(B_{ij}).$$

\begin{claim}\label{cl8}
For each $A_{ii},B_{ii} \in \mathcal{A}_{ii}$ such that $1\leq i \leq 2$, we have 
$$\Phi(A_{ii}+B_{ii})=\Phi(A_{ii})+\Phi(B_{ii}).$$
\end{claim}
We show that 
$$T=\Phi(A_{ii}+B_{ii})-\Phi(A_{ii})-\Phi(B_{ii})=0.$$
We can write
\begin{eqnarray*}
&&\Phi(A_{ii}+B_{ii})\diamond P_{j}+(A_{ii}+B_{ii})\diamond \Phi(P_{j})=\Phi((A_{ii}+B_{ii})\diamond P_{j})\\
&&=\Phi(A_{ii}\diamond P_{j})+\Phi(B_{ii}\diamond P_{j})\\
&&\Phi(A_{ii})\diamond P_{j}+A_{ii}\diamond \Phi(P_{j})+\Phi(B_{ii})\diamond P_{j}+B_{ii}\diamond \Phi(P_{j})\\
&&=(\Phi(A_{ii})+\Phi(B_{ii}))\diamond P_{j}+(A_{ii}+B_{ii})\diamond \Phi(P_{j}).
\end{eqnarray*}
So, we have
$$T\diamond P_{j}=0.$$
Therefore, we obtain $T_{ij}=T_{ji}=T_{jj}=0$.\\
On the other hand, for every $X_{ij}\in\mathcal{A}_{ij}$,  we have
\begin{eqnarray*}
&&\Phi(A_{ii}+B_{ii})\diamond X_{ij}+(A_{ii}+B_{ii})\diamond \Phi(X_{ij})=\Phi((A_{ii}+B_{ii})\diamond X_{ij})\\
&&=\Phi(A_{ii}\diamond X_{ij})+\Phi(B_{ii}\diamond X_{ij})=\Phi(A_{ii})\diamond X_{ij}+A_{ii}\diamond \Phi(X_{ij})\\
&&+\Phi(B_{ii})\diamond X_{ij}+B_{ii}\diamond \Phi(X_{ij})\\
&&=(\Phi(A_{ii})+\Phi(B_{ii}))\diamond X_{ij}+(A_{ii}+B_{ii})\diamond \Phi(X_{ij}).
\end{eqnarray*}
So, $$(\Phi(A_{ii}+B_{ii})-\Phi(A_{ii})-\Phi(B_{ii}))\diamond X_{ij}=0.$$
It follows that $T\diamond X_{ij}=0$  or $T_{ii}X_{ij}=0$. By knowing that $\mathcal{A}$ is prime, we have $T_{ii}=0$.\\
Hence, the additivity of $\Phi$ comes from the above claims.\\

In the rest of this paper, we show that $\Phi$ is $\ast$-derivation. 

\begin{claim}\label{cl9}
$\Phi$ preserves star.
\end{claim}
Since $\Phi(I)=0$ then we can write
$$\Phi\left(I\diamond A\right)=I\diamond \Phi(A).$$
Then 
$$\Phi(A-A^{*})=\Phi(A)-\Phi(A)^{*}.$$
So, we showed that $\Phi$ preserves star.

\begin{claim}\label{cl10}
we prove that $\Phi$ is derivation.
\end{claim}
For every $A,B\in\mathcal{A}$ we have
\begin{eqnarray*}
\Phi(AB-B^{*}A^{*})&=&\Phi(A^{*}\diamond B)\\
&=& \Phi(A^{*})\diamond B+A^{*}\diamond \Phi(B)\\
&=&\Phi(A^{*})^{*}B-\Phi(B)^{*}A^{*}-B^{*}\Phi(A^{*})+A\Phi(B).
\end{eqnarray*}
On the other hand, since $\Phi$ preserves star, we have
\begin{equation}\label{eq1cl10}
\Phi(AB-B^{*}A^{*})=\Phi(A)B+A\Phi(B)-B^{*}\Phi(A^{*})-\Phi(B)^{*}A^{*}.
\end{equation}
So, from (\ref{eq1cl10}), we have
\begin{eqnarray*}
&&\Phi(i(AB+B^{*}A^{*})=\Phi(A(iB)-(iB)^{*}A^{*})\\
&&=\Phi(A)(iB)+A\Phi(iB)-(iB)^{*}\Phi(A^{*})-\Phi(iB)^{*}A^{*}.
\end{eqnarray*}
Therefore, from claim \ref{cl3} we have
\begin{equation}\label{eq2cl10}
\Phi(AB+B^{*}A^{*})=\Phi(A)B+A\Phi(B)-B^{*}\Phi(A^{*})-\Phi(B^{*})A^{*}.
\end{equation}
By adding equations (\ref{eq1cl10}) and (\ref{eq2cl10}), we have
$$\Phi(AB)=\Phi(A)B+A\Phi(B).$$
This completes the proof.\\

\end{document}